\documentclass{amsart}
\usepackage{amsfonts}
\usepackage{amsmath}
\usepackage{amssymb}
\usepackage[utf8]{inputenc}
\usepackage{mathrsfs}
\newtheorem{theorem}{Theorem}
\newtheorem{corollary}{Corollary}
\newtheorem{lemma}{Lemma}
\newtheorem{proposition}{Proposition}

\newtheorem{definition}{Definition}
\newtheorem{remark}{Remark}

\begin{document}

\title[Graded Identities and Isomorphisms on Algebras of Upper Block-Triangular Matrices]{Graded Identities and Isomorphisms on Algebras of Upper Block-Triangular Matrices}

\author{Alex Ramos}
\address{Unidade Acadêmica de Matemática, Universidade Federal de Campina Grande, Campina Grande, PB, 58429-970, Brazil}
\email{almat100@hotmail.com}

\author{Diogo Diniz}
\address{Unidade Acadêmica de Matemática, Universidade Federal de Campina Grande, Campina Grande, PB, 58429-970, Brazil}
\email{diogo@mat.ufcg.edu.br}



\begin{abstract}
Let $G$ be an abelian group and $\mathbb{K}$ an algebraically closed field of characteristic zero. A. Valenti and M. Zaicev described the $G$-gradings on upper block-triangular matrix algebras provided that $G$ is finite. We prove that their result holds for any abelian group $G$: any grading is isomorphic to the tensor product $A\otimes B$ of an elementary grading $A$ on an upper block-triangular matrix algebra and a division grading $B$ on a matrix algebra. We then consider the question of whether graded identities $A\otimes B$, where $B$ is an algebra with a division grading, determine $A\otimes B$ up to graded isomorphism. In our main result, Theorem \ref{main}, we reduce this question to the case of elementary gradings on upper block-triangular matrix algebras which was previously studied by O. M. Di Vincenzo and E. Spinelli.
\end{abstract}
\maketitle

\section{Introduction}

Algebras of upper block-triangular matrices are important examples of p.i.-algebras. The verbal ideal of ordinary identities of an algebra upper block-triangular matrices is the product of the ideals of identities of each block (see \cite{GZ2}). Minimal varieties of finite basic rank were classified by A. Giambruno and M. Zaicev in terms of such algebras. The general result for arbitrary minimal varieties is in terms of of mimimal superalgebras, which may be realized as subalgebras of upper block-triangular matrix algebras with a $\mathbb{Z}_2$-grading (see \cite[Chapter 8]{GZ}). In this context it is of interest to classify the gradings by a group $G$ on an algebra of upper block-triangular matrices over a field $\mathbb{K}$ and to study its graded identities. For a finite abelian group and an algebraically closed field of characteristic zero the gradings were classified in \cite{VZ}. In Section \ref{grut} we prove, as a corollary of the result of A. Valenti and M. Zaicev, that given an algebra $UT(d_1,\dots, d_m)$ of upper block-triangular matrices with a grading by an abelian group $G$ there exist an algebra $UT(p_1,\dots, p_m)$ with an elementary grading and a matrix algebra $M_t(\mathbb{K})$ with a division grading such that $UT(d_1,\dots, d_m)$ is isomorphic as a graded algebra to $UT(p_1,\dots, p_m)\otimes M_t(\mathbb{K})$, here $\mathbb{K}$ is an algebraically closed field of characteristic zero. It was conjectured in \cite{VZ2} that this holds for any group $G$ and any field $\mathbb{K}$. This conjecture was recently solved in \cite{Y} with a minor restriction on the characteristic of the field. 

In Section \ref{id} we consider the following question: Let $R=A\otimes B$ and $R^{\prime}=A^{\prime}\otimes B^{\prime}$, where $A$ and $A^{\prime}$ are algebras of upper block-triangular matrices with elementary gradings and $B$ and $B^{\prime}$ are algebras with division gradings. If $R$ and $R^{\prime}$ satisfy the same graded identities then are $R$ and $R^{\prime}$ isomorphic as graded algebras? If $A$ is a matrix algebra then $R$ and $R^{\prime}$ are graded simple algebras. The graded Wedderburn-Artin theorem implies that any graded simple algebra that satisfies the d.c.c. on graded left ideals is the tensor product of a matrix algebra with an elementary grading and an algebra with a division grading. In \cite{AH}, \cite{PZ} it was proved that graded simple algebras are determined up to isomorphism by their graded identities. In \cite{DS} the authors consider algebras of upper block-triangular matrices with elementary gradings, it is proved in this paper that the answer is positive under suitable conditions. In our main result we reduce the previous question to the case of elementary gradings. As a consequence of the results of Di Vincenzo and Spinelli we conclude, for example, that the algebras $R$ and $R^{\prime}$ are isomorphic as graded algebras if the number of blocks in $A$ is $1$ or $2$. 

The main results of Section \ref{id} are that if $Id_G(R)=Id_G(R^{\prime})$ then $\mathrm{supp}\ B=\mathrm{supp}\ B^{\prime}$ and $Id_{H}(B)=Id_{H}(B^{\prime})$, where $H=\mathrm{supp}\ B$ (see Corollary \ref{supp} and Proposition \ref{div}). Now we assume that $H$ is a finitelly generated subgroup of $G$, since the field is algebraically closed this implies that $B$ is isomorphic to $B^{\prime}$ as a graded algebra. It was proved by D. S. Passmann in \cite{P} that a twisted group algebra of a solvable group over a field of characteristic zero is semisimple. We use this result to conclude that if a graded division algebra satisfies an ordinary polynomial identity then it satisfies the same ordinary identities as a matrix algebra (see Corollary \ref{iddiv}). As a consequence of this result we prove that if $B$ and $B^{\prime}$ satisfy an ordinary polynomial identity then the coarsenings $^{\alpha} R$ and $^{\alpha} R^{\prime}$, induced by the canonical quotient map $\alpha:G\rightarrow G/H$, satisfy the same identities as suitable algebras $U$ and $U^{\prime}$ (described in Theorem \ref{main}), respectively, of upper block-triangular matrices with elementary gradings (see Lemma \ref{elemc}). In our main result, Theorem \ref{main}, we prove that $Id_{G/H}(U)=Id_{G/H}(U^{\prime})$ and that if $U$ and $U^{\prime}$ are isomorphic as $G/H$-graded algebras then $R$ and $R^{\prime}$ are isomorphic as $G$-graded algebras.  

\section{Preliminaries}

Let $\mathbb{K}$ be a field and let $G$ be a group with identity element $e$. A grading by the group $G$ on a $\mathbb{K}$-algebra $A$ is a vector space decomposition $A=\oplus_{g\in G} A_g$ such that $A_gA_h\subseteq A_{gh}$ for every $g,h \in G$. An element $a$ of $A$ is homogeneous in the grading if $a\in A_g$ for some $g\in G$, if $a\neq 0$ we say that $g$ is the degree of $a$ in this grading and denote it $\mathrm{deg}_G\ a$. The subspace $A_e$ is the neutral component of $A$, clearly it is a subalgebra of $A$. The support of the grading is the set \[\mathrm{supp}\ A =\{g\in G\mid A_g\neq 0\}.\] 

The grading on $A$ is a division grading if $A$ has a unit and every non-zero homogeneous element is invertible. If $\mathbb{K}$ is an algebraically closed field then any graded division algebra is isomorphic to a twisted group algebra $\mathbb{K}^{\alpha}G$, where $\alpha$ is a $2$-cocycle (see \cite[Theorem 2.13]{EK}). The map $\beta_{\alpha}:G\times G \rightarrow \mathbb{K}^{\times}$ given by $$\beta_{\alpha}(g,h)=\frac{\alpha(g,h)}{\alpha(h,g)}$$ is an alternating bicharacter. If $G$ is a finitely generated abelian group it is easy to verify that two twisted group algebras $\mathbb{K}^{\alpha}G$ and $\mathbb{K}^{\alpha^{\prime}}G$ with the same associated bicharacter are isomorphic. Moreover, it is clear that two graded division algebras with the same graded identities have the same associated bicharacter. This proves the following remark.

\begin{remark}\label{idisgdiv}
	Let $G$ be a finitely generated abelian group and let $\mathbb{K}$ be an algebraically closed field. If  $D_1$ and $D_2$ be two $\mathbb{K}$-algebras with division $G$-gradings and $Id_G(D_1)=Id_{G}(D_2)$ then $D_1$ is isomorphic to $D_2$ as a graded algebra.
\end{remark} 

Given a homomorphism $\alpha:G\rightarrow H$ of groups and a $G$-grading $A=\oplus_{g\in G}A_g$ on the algebra $A$ it follows that the decomposition $A=\oplus_{h\in H} A_h$, where $A_h=\oplus_{g \in \alpha^{-1}(h)} A_g$, is an $H$-grading on $A$. We refer to this as the grading induced by $\alpha$ and denote it $^{\alpha} A$

Let $A$ and $B$ be algebras graded by the group $G$. A map $f:A\rightarrow B$ is a homomorphism of graded algebras if it is a homomorphism of algebras such that $f(A_g)\subseteq B_g$ for every $g\in G$. If $f$ is also bijective then we say it is an isomorphism of graded algebras.

Given an $m$-tuple of positive integers $(d_1,\dots, d_m)$ we denote $UT(d_1,\dots, d_m)$ the subalgebra of $M_n(\mathbb{K})$, where $n=d_1+\cdots +d_m$, of the matrices of the form \[\left(\begin{array}{ccc}
A_{1 1}&\cdots & A_{1 m}\\
\vdots & \ddots & \vdots \\
0&\cdots & A_{m m}
\end{array}\right),\] where $A_{i j}$ is a block of size $d_i\times d_j$. We refer to $UT(d_1,\dots, d_m)$ as an algebra of upper block-triangular matrices.  A grading by a group $G$ on $A=UT(d_1,\dots, d_m)$ is elementary if every elementary matrix in $A$ is homogeneous in the grading. In this case there exist an $n$-tuple $(g_1,\dots, g_n)\in G^{n}$ such that the elementary matrix $e_{ij}\in A$ is homogeneous of degree $g_ig_j^{-1}$. Moreover if $B$ is any $G$-graded algebra then $A\otimes B$ has a $G$-grading such that $(A\otimes B)_g$ is the subspace spanned by the elements $e_{ij}\otimes b$ such that $b$ is a homogeneous element in $B$ and $g_i(\mathrm{deg}_G\ b)g_j^{-1}=g$.

The free $G$-graded algebra is the free algebra $\mathbb{K}\langle X_G \rangle$, freely generated by the set $X_G$ which is the disjoint union of the sets $X_g=\{x_{1}^g,x_{2}^g,\cdots\}$, $g\in G$, with the grading in which the variables in $X_g$ are homogeneous of degree $g$. In order to simplify the notation we may omit the neutral element $e$ in the notation, thus we denote by $X$ the set $X_e$ and by $x_i$ the indeterminate $x_{i}^e$ in $X$. Let $f(x_{1}^{g_1},\dots, x_{n}^{g_n})$ be a polynomial in $\mathbb{K}\langle X_G \rangle$ and $A=\oplus_{g\in G} A_g$ be a $G$-graded algebra. An admissible substitution for $f$ is an $n$-tuple $(a_1,\dots, a_n)$ of elements of $A$ such that $a_i\in A_{g_i}$ for $i=1,\dots, n$. We say that $f$ is a graded polynomial identity for $A$ if $f(a_1,\dots, a_n)=0$ whenever $(a_1,\dots, a_n)$ is an admissible substitution for $f$.  Note that $f\in K\langle X_G \rangle$ is a graded polynomial identity for $A$ if and only if $f$ lies in the kernel of every graded homomorphism $K\langle X_G \rangle\rightarrow A$. The set of graded polynomial identities for $A$ is denoted $Id_G(A)$. This is an ideal of $K\langle X_G \rangle$ invariant for graded endomorphisms of  $K\langle X_G \rangle$. The subalgebra $\mathbb{K}\langle X \rangle$ of $\mathbb{K}\langle X_G \rangle$ is the ordinary free algebra. A polynomial $f(x_1,\dots, x_n)\in \mathbb{K}\langle X \rangle$ is a polynomial identity for an algebra $B$ if $f(b_1,\dots, b_n)=0$ for every $b_1,\dots, b_n \in B$. We denote $Id(B)$ the ideal of identities of $B$. 

\section{Gradings on Upper Block-Triangular Matrix Algebras}\label{grut}

In this section we present the classification of the group gradings on upper block-triangular matrix algebras up to isomorphism. The base field is algebraically closed of characteristic zero and the group is abelian. The main result of \cite{VZ} is the following theorem.

\begin{theorem}\cite[Theorem 3.2]{VZ}\label{VZ}
	Let $G$ be a finite abelian group and $UT(d_1,\dots, d_m)$ an upper block-triangular matrix algebra over an algebraically closed field $\mathbb{K}$ of characteristic zero. Then there exist a decomposition $d_1 = p_1t, \dots, d_m=p_mt$, a subgroup $H\subset G$, and an $n$-tuple $(g_1, \dots, g_n)\in G^n$ , where $n = p_1 + \cdots + p_m$, such that $UT(d_1,\dots, d_m)$ is isomorphic to $UT(p_1,\cdots , p_m)\otimes M_t(\mathbb{K})$ as a $G$-graded algebra where $M_t(\mathbb{K})$ is an $H$-graded algebra with a “fine” grading with support $H$ and $UT(p_1,\cdots , p_m)$ has an elementary grading defined by $(g_1, \dots, g_n)$.
\end{theorem}

We prove, as a corollary, that the above theorem holds for arbitrary abelian groups.

\begin{lemma}\label{grading}
	Let $\varphi:G\rightarrow H$ be a homomorphism of groups and let $S$ be a subset of $G$ such that the restriction of $\varphi$ to $S\cup (S\cdot S)$ is injective. Let $A$ be a $G$-graded algebra such that $(\mathrm{supp}\ A )\subseteq S$ and let $B$ be an $H$-graded algebra such that $(\mathrm{supp}\ B)\subseteq \varphi(S)$.
	\begin{enumerate}
	 \item[(i)] If $B_s= B_{\varphi(s)}$ for any $s\in  S$ and $B_g=0$ for any $g\in G\setminus S$ then $B=\oplus_{g\in G}B_g$ is a $G$-grading on $B$. A subspace of $B$ is homogeneous in this  $G$-grading if and only if it is homogeneous in the $H$-grading.
	 \item[(ii)] An isomorphism of algebras $f:A\rightarrow B$ is an isomorphism of $G$-graded algebras from $A$ to $B$ with this $G$-grading if and only if $f$ is an isomorphism of $H$-graded algebras form $^{\varphi}A$ to $B$ with its $H$-grading.
	\end{enumerate}
\end{lemma}
\textit{Proof.}
If $g\in G\setminus S$ or $h \in G\setminus S$ then $B_gB_h=0\subseteq B_{gh}$. It remains to prove that $B_sB_t\subseteq B_{st}$ for every $s,t\in S$. In this case $B_s=B_{\varphi(s)}$ and $B_t=B_{\varphi(t)}$, therefore
\begin{equation}\label{st}
B_sB_t=B_{\varphi(s)}B_{\varphi(s)}\subseteq B_{\varphi(st)}.
\end{equation}
If $B_{\varphi(st)}=0$ then $B_sB_t=0$ and therefore $B_sB_t\subseteq B_{st}$. Now we assume that $B_{\varphi(st)}\neq 0$. In this case $\varphi(st)\in \mathrm{supp}\ B \subseteq \varphi(S)$. Therefore there exists an $s^{\prime}\in S$ such that $\varphi(st)=\varphi(s^{\prime})$. Since the restriction $\varphi$ to $S\cup(S\cdot S)$ is injective we conclude that $st=s^{\prime}$. Hence $$B_{\varphi(st)}=B_{\varphi(s^{\prime})}=B_{s^{\prime}}=B_{st}.$$ Then it follows from (\ref{st}) that $B_sB_t\subseteq B_{st}$. An element of $b\in B$ is homogeneous in the $G$-grading if and only if $b$ is homogeneous in the $H$-grading, this implies that a subspace $V$ is homogeneous in the $G$-grading of $B$ if and only if it is homogeneous in the $H$-grading. This proves [(i)]. Note that if $b$ is a homogeneous element of $B$ we have $\mathrm{deg}_H\ b=\varphi(\mathrm{deg}_G\ b)$. Since $(\mathrm{supp}\ A )\subseteq S$ and the restriction of $\varphi$ to $S$ is injective we conclude that an element $a\in A$ is homogeneous in the $G$-grading if and only if $a$ is homogeneous in the $H$-grading $^{\varphi}A$ and that $\mathrm{deg}_H\ a=\varphi(\mathrm{deg}_G\ a)$. This implies [(ii)].
\hfill $\Box$

Next we prove that the above theorem holds for arbitrary abelian groups. We remark that in \cite{Y} it was proved that this result holds in general provided that the base field $\mathbb{K}$ is of characteristic zero or that $\mathrm{char}\ \mathbb{K}$ is strictly greater than $\mathrm{dim}\ UT(d_1,\dots, d_m)$.

\begin{corollary}\label{garb}
	The previous theorem holds for abelian groups.
\end{corollary}
\textit{Proof.}
 Let $G$ be an abelian group and $\mathbb{K}$ an algebraically closed field of characteristic zero. Let $A=UT(d_1,\dots, d_m)$ be an upper block triangular matrix algebra over $\mathbb{K}$ with a $G$-grading. We assume without loss of generality that $G$ is generated by $S=\mathrm{supp}\ A$. Since $A$ has finite dimension the support $S$ is finite, therefore $G$ is isomorphic to a direct product of cyclic groups $\mathbb{Z}^p\times \mathbb{Z}_{n_1}\times \cdots \times \mathbb{Z}_{n_q}$, where $p,q\geq 0$ and $n_1,\dots, n_q>0$ are integers. Hence we may assume without loss of generality that $G=\mathbb{Z}^p\times \mathbb{Z}_{n_1}\times \cdots \times \mathbb{Z}_{n_q}$. Given an integer $m>0$ let $H=(\mathbb{Z}_m)^p \times \mathbb{Z}_{n_1}\times\cdots \times \mathbb{Z}_{n_q}$. Denote $\varphi:G\rightarrow H$ the homomorphism given by $$\varphi(z_1,\dots, z_p,w_1,\dots, w_q)=(\overline{z_1},\dots, \overline{z_p},w_1,\dots, w_q),$$ where $z_1,\dots, z_p\in \mathbb{Z}$, $w_i\in \mathbb{Z}_{n_i}$, for $i=1,\dots, q$ and $z\mapsto \overline{z}$ is the canonical homomorphism $\mathbb{Z}\rightarrow \mathbb{Z}_m$. Given any finite subset $T$ of $G$ we may choose $m$ large enough so that the restriction of $\varphi$ to $T$ is injective. Since the support $S$ is finite this implies that we may choose $m$ large enough so that the restriction of $\varphi$ to $T=S\cup (S\cdot S)$ is an injective map. Let $^{\varphi}A$ be the $H$-grading induced by $\varphi$. Theorem \ref{VZ} implies that the coarsening $^{\varphi}A$ is isomorphic as an $H$-graded algebra, to $B=UT(p_1,\cdots , p_m)\otimes M_t(\mathbb{K})$, where $UT(p_1,\cdots , p_m)$ has an elementary grading and $M_t(\mathbb{K})$ has a fine grading. Denote $B=\oplus_{h\in H} B_h$ the $H$-grading on $B$. Let $B_s= B_{\varphi(s)}$ for any $s\in  S$ and $B_g=0$ for any $g\in G\setminus S$. Lemma \ref{grading} implies that $B=\oplus_{g\in G} B_g$ is a $G$-grading on $B$. Moreover any isomorphism of $H$-graded algebras from $^{\varphi}A$ to $B$ is also an isomorphism of $G$-graded algebras from $A$ to $B$ with this $G$-grading, hence $A$ with its $G$-grading is isomorphic as a graded algebra to $B$ with this $G$-grading. A subspace of $B$ is homogeneous in this  $G$-grading if and only if it is homogeneous in the $H$-grading, therefore we obtain a fine $G$-grading on $M_t(\mathbb{K})$ and an elementary $G$-grading on $UT(p_1,\dots, p_m)$ such that $B$ with the $G$-grading is $UT(p_1,\dots, p_m)\otimes M_t(\mathbb{K})$ with the canonical $G$-grading on the tensor product.
\hfill $\Box$

In \cite{ACD} the isomorphism classes of upper block-triangular matrix algebras graded by an abelian group were described. 

\begin{corollary}\cite[Corollary 4]{ACD}\label{appl}
	Let $G$ be an abelian group and let $\mathbb{K}$ be an arbitrary field. Let $R$ and $R^{\prime}$ denote the algebras $ UT(p_1,\dots , p_m)\otimes_{\mathbb{K}} B$ and $ UT(p_1^{\prime},\dots , p_m^{\prime})\otimes_{\mathbb{K}}B^{\prime}$ respectively, where $B$, $B^{\prime}$ are graded division algebras and $UT(p_1,\dots , p_m)$, $UT(p_1^{\prime},\dots , p_m^{\prime})$ have elementary gradings defined by the tuples $\textup{\textbf{g}}$, $\textup{\textbf{g}}^{\prime}$ of elements of $G$, respectively. The algebras $R$ and $R^{\prime}$ are isomorphic if and only if $B\cong B^{\prime}$, $(p_1,\dots , p_m)=(p_1^{\prime},\dots , p_m^{\prime})$ and there exist a $g\in G$, $h_1,\dots, h_n\in \mathrm{supp}\ B$ and $\sigma\in S_{p_1}\times \dots \times S_{p_m}$ such that $g_i^{\prime}=g_{\sigma(i)}h_{\sigma(i)}g$ for $i=1,\dots,n$.
\end{corollary}

We remark that in \cite{ACD} the previous corollary was sated for $B=M_t(\mathbb{K})$ and $B^{\prime}=M_{t^{\prime}}(\mathbb{K})$, this however is not used in the proof.

Let $R=A\otimes B$, where $A=UT(d_1,\dots , d_m)$ is an upper block-triangular matrix algebra with an elementary grading and $B$ is a graded division algebra with support $H$. Let ${\bf g}=(g_1,\dots, g_n)$ be a tuple of elements of $G$ that determines the elementary grading in $A$. Now let $g_1^{\prime},\dots, g_n^{\prime}$ representatives of the cosets $g_1H,\dots, g_nH$, respectively, such that, $g_iH=g_i^{\prime}H$ for $i=1,\dots, n$ and $g_i^{\prime}=g_j^{\prime}$ whenever $g_i^{\prime}H=g_j^{\prime}H$. Let $A^{\prime}$ denote the algebra $UT(d_1,\dots,d_m)$ with the elementary grading induced by ${\bf g}^{\prime}=(g_1^{\prime},\dots, g_n^{\prime})$. It follows from Corollary \ref{appl} that the algebras $R$ and $R^{\prime}=A^{\prime}\otimes B$ are isomorphic. Since $(g_i^{\prime})^{-1}g_j^{\prime}\in H$ if and only if $g_i^{\prime}=g_j^{\prime}$ we conclude that $R^{\prime}_e=A^{\prime}_e\otimes 1\cong A_e^{\prime}$.

\begin{remark}\label{ecomp}
	Let $R=A\otimes B$, where $A=UT(d_1,\dots , d_m)$ is an upper block-triangular matrix algebra with an elementary grading and $B$ is a graded division algebra. The algebra $UT(d_1,\dots , d_m)$ may be equipped with an elementary grading $A^{\prime}$  such that $R$ is isomorphic to $R^{\prime}$, where $R^{\prime}=A^{\prime}\otimes B$ and $R^{\prime}_e=A_e^{\prime}\otimes 1$.
\end{remark}

\section{Graded Identities for Algebras of Upper Block-Triangular Matrices}\label{id}

In this section we present the main results of the paper. Next we define the Capelli polynomials, these polynomials will play an important role in this section.

\begin{definition}
	The Capelli polynomial of rank $t$ is the polynomial \[Cap_t(x_1,\dots, x_t, y_1,\dots, y_{t+1})=\sum_{\sigma \in S_t}(\mathrm{sgn}\ \sigma) y_1x_{\sigma(1)}y_2x_{\sigma(2)}\cdots y_t x_{\sigma(t)}y_{t+1},\] in $\mathbb{K}\langle X \rangle$.
\end{definition}

\begin{lemma}\cite[Lemma 9.1.4]{GZ}\label{cap}
	If $d=d_1^2+\cdots d_m^2$, the algebra $UT(d_1,\dots, d_m)$ satisfies the Capelli identity $\mathrm{Cap}_{d+m}\equiv 0$ but does not satisfy $\mathrm{Cap}_{d+m-1}\equiv 0$.
\end{lemma}

\begin{remark}\label{jacob}
	As noted in \cite[Lemma 3.2]{DS} it follows from the proof of Lemma 9.1.4 in \cite{GZ} that the result of any substitution in $\mathrm{Cap}_{d+m-1}$ lies in $J(A)^{m-1}$, where $A=UT(d_1,\dots, d_m)$ and $J(A)$ denotes the Jacobson radical of $A$. 
\end{remark}

\begin{proposition}\label{decomp}
Let $A=UT(d_1,\dots, d_m)$ denote an upper-block triangular matrix algebra with an elementary grading induced by ${\bf g}=(g_1,\dots, g_n)$, where $n=d_1+\cdots + d_m$. Let $r$ denote the number of distinct elements that appear in ${\bf g}$. Then the neutral component $A_e$ is a direct sum of $r$ ideals, each isomorphic to an upper-block triangular matrix algebra of length $\leq m$. 
\end{proposition}
\textit{Proof.}
Let $I_1=\{1,\dots,d_1\}$ and $$I_{t}=\{d_1+\cdots + d_{t-1}+1,\dots,d_1+\cdots +d_{t-1}+ d_t\},$$ for $t=2,\dots,m$. The sets $I_1,\dots, I_m$ form a partition of $N=\{1,\dots,n\}$. A basis for $A$ consists of the matrix units $e_{ij}$ where $i\in I_s$, $j\in I_t$ and $s\leq t$. Denote $h_1,\dots, h_r$ the elements of $G$ that appear in the $n$-tuple ${\bf g}$, moreover let $J_k=\{i\in N \mid g_i=h_k\}$. It is clear that a basis for $A_e$ consists of the matrix units $e_{ij}$ in the basis for $A$ such that $i,j\in J_k$ for some $k$. Therefore $$A_e=B_1\oplus\cdots \oplus B_r,$$ where $B_k$ is the subspace of $A_e$ generated by the matrix units $e_{ij}$ in the basis for $A$ such that $i,j\in J_k$. It is easy to verify that $B_k$ is an ideal in $A_e$. A basis for $B_k$ is the set of matrix units $e_{ij}$ where $i\in I_s\cap J_k$, $j\in I_t \cap J_k$ and $s\leq t$, therefore $B_k$ is isomorphic an algebra of upper block-triangular matrices of length $\leq m$. 
\hfill $\Box$

The following remarks will be used in the proof of Lemma \ref{prev}.


\begin{remark}\label{elem}
	Let $a_1,\dots, a_t, v_1,\dots, v_{t-1}$ be matrix units in a matrix algebra and assume that $v_1,\dots, v_{t-1}$ are pairwise distinct. If  $$a_1v_1\cdots a_{t-1}v_{t-1}a_t\neq 0$$ then for every permutation $\tau \in S_{t-1}$ different from $Id$ we have $$a_1v_{\tau(1)}\cdots a_{t-1}v_{\tau(t-1)}a_t=0.$$
\end{remark}

\begin{lemma}\label{prev}
Let $A$ be an upper-block triangular matrix algebra with an elementary $G$-grading. Let $t$ be the natural number such that $\mathrm{Cap}_{t-1}$ is not an identity for $A_e$ and $\mathrm{Cap}_t$ is an identity for $A_e$. Then  $\mathrm{Cap}_{t-1}x_{t}^g\notin Id_G(A)$  and $$\mathrm{Cap}_t(x_1,\dots,x_{t-1},x_t^g,y_1,\dots,y_{t+1})\in Id_G(A)$$ if and only if $g=e$.
\end{lemma}
\textit{Proof.}
It is clear, from the choice of $t$, that if $g=e$ then $\mathrm{Cap}_t$ is an identity for $A$ and $\mathrm{Cap}_{t-1}x_{t}^g\notin Id_G(A)$. To prove the converse we assume that $\mathrm{Cap}_{t-1}x_{t}^g$ is not an identity for $A$ and that $g\neq e$ and prove that $\mathrm{Cap}_t(x_1,\dots,x_{t-1},x_t^g,y_1,\dots,y_{t+1})$ is not an identity for $A$. Let $A_e=B_1\oplus\cdots \oplus B_r$ be the decomposition in Proposition \ref{decomp}. Let $a_1,\dots, a_{t-1},a_t,v_1,\dots,v_{t-1}, v_t$ be an admissible substitution by matrix units such  that $$\mathrm{Cap}_{t-1}(v_1,\dots, v_{t-1},a_1,\dots,a_t)v_t\neq 0.$$ This implies that elements $v_1,\dots, v_{t-1}$ are pairwise distinct, moreover there exists $i$ such that $v_1,\dots, v_{t-1},a_1,\dots,a_t\in B_i$. We may reorder the elements $v_1,\dots,v_{t-1}$, if necessary, and assume that $$c:=a_1v_1\cdots a_{t-1}v_{t-1}a_tv_t\neq 0.$$  Clearly there exists $a_{t+1}\in A_e$ such that $ca_{t+1}\neq 0$. Since $\mathrm{deg}_G(v_{t})\neq e$ and $\mathrm{deg}_G(v_{i})=e$, for $i=1,\dots, t-1$, we conclude that the elements $v_1,\dots, v_t$ are pairwise distinct. Remark \ref{elem} implies that 
\begin{equation}\label{e}
a_1v_{\sigma(1)}a_2\cdots a_{t-1}v_{\sigma(t)}a_{t+1}= 0
\end{equation}
 whenever $\sigma \neq 1$. Therefore we conclude that $\mathrm{Cap}_t(v_1,\dots,v_{t},a_1,\dots,a_{t+1})\neq 0$.
\hfill $\Box$

\begin{corollary}\label{supp}
Let $R=A\otimes B$, where $A$ is an upper block-triangular matrix algebra with an elementary grading and $B$ is a graded division algebra. Let $t$ be the natural number such that $\mathrm{Cap}_{t-1}$ is not an identity for $R_e$ and $\mathrm{Cap}_t$ is an identity for $R_e$. The element $g$ lies in $\mathrm{supp}\ B$ if and only if $\mathrm{Cap}_{t-1}x_{t}^g \notin Id_G(R)$ and the polynomial $\mathrm{Cap}_t(x_1,\dots,x_{t-1},x_t^g,y_1,\dots,y_{t+1})$ lies in $Id_G(R)$.
\end{corollary}

\textit{Proof.}
 Remark \ref{ecomp} implies that we may assume without loss of generality that $R_e=A_e\otimes 1$. 
 
 ($\Rightarrow$)
 
 Let $g$ be an element of $\mathrm{supp}\ B$. Since $\mathrm{Cap}_{t-1}$ is not an identity for $R_e$ exist $c_1,\dots, c_{t-1}, d_1,\dots, d_t\in A_e$ such that $$\mathrm{Cap}_{t-1}(c_1,\dots, c_{t-1},d_1,\dots, d_t)\neq 0.$$ Let $b$ be a non-zero element in $B_g$. We have $$\left(\mathrm{Cap}_{t-1}(c_1\otimes 1,\dots, c_{t-1}\otimes 1,d_1\otimes 1,\dots, d_t\otimes 1)\right)(1\otimes b) \neq 0.$$ Hence $\mathrm{Cap}_{t-1}x_{t}^g \notin Id_G(R)$. Since $1\otimes b$ is a homogeneous invertible element of degree $g$ in $R$ and $R_e=A_e\otimes 1$ we conclude that $R_g=\{a\otimes b \mid a \in A_e\}$. Therefore given a homomorphism of graded algebras $\varphi:F\langle X_G \rangle \rightarrow R$ there exist $v_1,\dots, v_{t-1}, v_t,w_1,\dots, w_{t+1}\in A_e$ such that $\varphi(x_i)=v_i\otimes 1$, $i=1,\dots, t-1$, $\varphi(x_t^{g})=v_t\otimes b$ and $\varphi(y_j)=w_j\otimes 1$, $j=1,\dots, t+1$. Hence 
 \[
 \begin{array}{cc}
 &\varphi\left(\mathrm{Cap}_t(x_1,\dots,x_{t-1},x_t^g,y_1,\dots,y_{t+1})\right)=\\&\mathrm{Cap}_t(v_1,\dots,v_{t-1},v_t,w_1,\dots,w_{t+1})\otimes b=0. 
 \end{array}
 \]
 Therefore $\mathrm{Cap}_t(x_1,\dots,x_{t-1},x_t^g,y_1,\dots,y_{t+1})$ lies in $Id_G(R)$.
 
 ($\Leftarrow$)
 
Since $\mathrm{Cap}_{t-1}x_{t}^g\notin Id_G(R)$ there exists a homomorphism of graded algebras $\varphi:F\langle X_G \rangle \rightarrow R$ such that $\varphi(\mathrm{Cap}_{t-1}x_{t}^g)\neq 0$. We write $\varphi(x_{t}^g)=\sum_{i=1}^n a_i\otimes b_i$, where $a_i\in A$ and $b_i\in B$ are homogeneous elements and 
\begin{equation}\label{deg}
(\mathrm{deg}_G a_i)\cdot (\mathrm{deg}_G b_i)=g,
\end{equation} for $i=1,\dots, n$. Since $R_e=A_e\otimes 1$ there exist $c_1,\dots, c_{t-1}, d_1,\dots, d_t\in A_e$ such that $$\varphi(\mathrm{Cap}_{t-1})=\mathrm{Cap}_{t-1}(c_1,\dots, c_{t-1},d_1,\dots, d_t)\otimes 1.$$ Note that $\varphi(\mathrm{Cap}_{t-1})\cdot(a_j\otimes b_j)\neq 0$ for some $j$. This implies that $$\mathrm{Cap}_{t-1}(c_1,\dots, c_{t-1},d_1,\dots, d_t)\cdot a_j\neq 0,$$ hence  the polynomial $\mathrm{Cap}_{t-1}x_{t}^{g_1}$, where $g_1=(\mathrm{deg}_G\ a_j)$, is not a graded identity for $A$. Given $v_1,\dots, v_{t-1},w_1,\dots, w_{t+1}\in A_e$ and $v_t\in A_{g_1}$ we have
\[
\begin{array}{ccl}
	0&=&\mathrm{Cap}_t(v_1\otimes 1,\dots,v_{t-1}\otimes 1,v_t\otimes b_j,w_1\otimes 1,\dots,w_{t+1}\otimes 1)\\ &=&\mathrm{Cap}_t(v_1,\dots,v_{t-1},v_t,w_1,\dots,w_{t+1})\otimes b_j.
\end{array}\] The element $b_j$ is a non-zero homogeneous element in $B$, hence it is invertible. This implies that $\mathrm{Cap}_t(v_1,\dots,v_{t-1},v_t,w_1,\dots,w_{t+1})=0$. Therefore the polynomial $\mathrm{Cap}_t(x_1,\dots,x_{t-1},x_t^{g_1},y_1,\dots,y_{t+1})$ is a graded identity for $A$. The previous lemma implies that $g_1=e$, therefore it follows from (\ref{deg}) that $g=(\mathrm{deg}_G b_i)\in \mathrm{supp}\ B$.
\hfill $\Box$

\begin{lemma}\label{1}
Let $A=UT(d_1,\dots, d_m)$ be an upper block-triangular matrix algebra and $t=d_1^{2}+\cdots d_m^2+m$. Let $f$ be a polynomial such that  $f\mathrm{Cap}_{t-1}$ is a multilinear polynomial. If the result of a substitution is a non-zero element in $A$ then the variables in $f$ are replaced by elements in the first block of $A$.
\end{lemma}

\textit{Proof.}
Since $f\mathrm{Cap}_{t-1}$ is a multilinear polynomial we only need to consider substitutions of the variables by matrix units. Let $\varphi:F\langle X \rangle \rightarrow A$ be a homomorphism such that $\varphi(x)$ is a matrix unit for every $x\in X$ and $\varphi(f\mathrm{Cap}_{t-1})\neq 0$. Remark \ref{jacob} implies that $\varphi(\mathrm{Cap}_{t-1})\in J(A)^{m-1}$. Therefore $\varphi(\mathrm{Cap}_{t-1})=e_{rs}$ with $1\leq r \leq d_1$. Since $J(A)^m=0$, if any variable in $f$ is replaced by an element in $J(A)$ we have $\varphi(f)\in J(A)$. Therefore $\varphi(f\mathrm{Cap}_{t-1})$ lies in $J(A)^m=0$, a contradiction. Hence the variables in $f$ are replaced by elements in the diagonal blocks of $A$. If variables are replaced by elements in distinct blocks then $\varphi(f)=0$. Thus the variables in $f$ are replaced by elements in a single diagonal block. Hence $\varphi(f)=e_{tu}$ with $d_1+\cdots +d_{i-1}+1\leq t,u \leq d_1+\cdots +d_i$. We have $$0\neq \varphi(f\mathrm{Cap}_{t-1})=e_{tu}e_{rs},$$ hence $u=r$. Therefore we conclude that $i=1$.
\hfill $\Box$

\begin{proposition}\label{fcap}
	Let $U=U_1\oplus\cdots \oplus U_r$, where $U_1,\dots, U_r$ are ideals of $U$ isomorphic to algebras of upper block-triangular matrices. Denote $t$ the least integer such that $\mathrm{Cap}_t\in Id(U)$. There exists a polynomial $f$ in $\mathbb{K}\langle X \rangle$ such that: 
	\begin{enumerate}
	\item[(i)] $f\mathrm{Cap}_{t-1}$ is a multilinear polynomial;
	\item[(ii)] $f\mathrm{Cap}_{t-1}$ is not an identity for $U$;
	\item[(iii)] If $\varphi:\mathbb{K}\langle X \rangle \rightarrow U$ is a homomorphism and $\varphi(f\mathrm{Cap}_{t-1})\neq 0$ then $\varphi(f)$ is in the center of $U$.
	\end{enumerate}	  
\end{proposition}
\textit{Proof.}
We may assume without loss of generality that have $U_i=UT(d_{1}^{(i)},\dots, d_{m_i}^{(i)})$, $i=1,\dots, r$. Denote $t_i=(d_1^{(i)})^2+\dots + (d_{m_i}^{(i)})^2+m_i$. Note that $t=\mathrm{max}\ \{t_1,\dots, t_r\}$. Denote $k=\mathrm{max}\ \{d_1^{(j_1)},\dots, d_1^{(j_s)}\}$, where $\{j_1,\dots,j_s\}=\{j\mid t_j=t\}$. Let $f$ be a central polynomial for $M_k(\mathbb{K})$ such that $f\mathrm{Cap}_{t-1}$ is a multilinear polynomial. It is clear that is not an identity for $U$. We claim that $(iii)$ holds. Let $\varphi:F\langle X \rangle \rightarrow U$ be a homomorphism such that $\varphi(f\mathrm{Cap}_{t-1})\neq 0$. Note that $\varphi(f\mathrm{Cap}_{t-1})$ is the result of the substitution $x\mapsto \varphi(x)$ in $f\mathrm{Cap}_{t-1}$. Since $\varphi(f\mathrm{Cap}_{t-1})\neq 0$ there exists an $i$ such that each indeterminate in $f\mathrm{Cap}_{t-1}$ is replaced by an element in $U_i$. Lemma \ref{1} implies that the variables in $f$ are replaced by elements in the first block of $U_i$. Recall that $f$ is an identity for $M_l(\mathbb{K})$ if $k>l$ (see \cite[pg. 150, Lemma 3.5]{DF}). Since $\varphi(\mathrm{Cap}_{t-1})\neq0$ we conclude that $t_i=t$. Therefore the choice of $k$ implies that $d_1^{(i)}=k$. Hence $\varphi(f)$ is a central element.
\hfill $\Box$

\begin{definition}
	Let $A$ be an upper block-triangular matrix algebra with an elementary grading and let $t$ be the least integer such that $\mathrm{Cap}_t\in Id(A_e)$. We say a multilinear polynomial $f$ is $A$-good if the polynomial $fCap_{t-1}$ is multilinear, $fCap_{t-1}$ is not a polynomial identity for $A_e$ and for every homomorphism $\varphi:\mathbb{K}\langle X \rangle \rightarrow A_e$  such that $\varphi(f\mathrm{Cap}_{t-1})\neq 0$ the element $\varphi(f)$ is in the center of $A_e$.
\end{definition}

\begin{remark}
	Proposition \ref{decomp} and  Proposition \ref{fcap} imply that for any upper block-triangular matrix algebra $A$ with an elementary grading there exists an $A$-good polynomial. 
\end{remark}

\begin{proposition}\label{div}
Let $R=A\otimes B$, where $A$ is an upper block-triangular matrix algebra with an elementary grading and $B$ is a graded division algebra. Let $t$ the least integer such that $\mathrm{Cap}_t\in Id(A_e)$ and let $f$ be an $A$-good polynomial. Let $S_1=\{y_1^{h_1},y_2^e,\dots, y_m^e\}, \dots, S_u=\{z_1^{h_u},z_2^e,\dots, z_m^e\}$ be $u$ pairwise disjoint sets of variables and let $$f^{h_1}=f(y_1^{h_1},y_2^e,\dots, y_m^e), \dots, f^{h_u}=f(z_1^{h_u},z_2^e,\dots, z_m^e).$$ A multilinear polynomial $g(x_1^{h_1},\dots, x_u^{h_u})\in F\langle X_H\rangle$, where $H=\mathrm{supp}\ B$, is a graded identity for $B$ if and only if the multilinear polynomial $g(f^{h_1},\dots, f^{h_{u}})\mathrm{Cap}_{t-1}$ is a graded identity for $R=A\otimes B$.
\end{proposition}
\textit{Proof.}

As a consequence of Remark \ref{ecomp} we may assume without loss of generality that $R_e=A_e\otimes1$. 

($\Rightarrow$)

For each $h\in H$ we fix a non-zero element $b_h \in B_h$. Since $R_h=R_e(1\otimes b_h)$ any element $a$ in $R_h$ is of the form $a_e\otimes b_h$, where $a_e\in A_e$.  Let $\varphi:F\langle X_G \rangle \rightarrow R$ be a homomorphism of graded algebras. Note that for each $l$ there exists an $f_l\in f(A_e,\dots, A_e)$ such that $\varphi(f^{h_l})=f_l\otimes b_{h_l}$. Moreover $\varphi(\mathrm{Cap}_{t-1})=d\otimes 1$ for some $d\in \mathrm{Cap}_{t-1}(A_e,\dots, A_e)$. The polynomial $f$ is an $A$-good polynomial, therefore either $f_l$ is a central element in $A_e$ or $f_ld=0$. Hence we conclude that $\varphi(g(f^{h_1},\dots, f^{h_{u}})\mathrm{Cap}_{t-1})=0$ if $f_ld=0$ for some $l$. It remains to consider the case where $f_1,\dots, f_u$ are central elements of $A_e$. In this case $\varphi(g(f^{h_1},\dots, f^{h_{u}}))=f_1\cdots f_u\otimes g(b_{h_1},\cdots, b_{h_u})=0$, because $g(x_1^{h_1},\dots, x_u^{h_u})\in Id_H(B)$. Therefore $g(f^{h_1},\dots, f^{h_{u}})\mathrm{Cap}_{t-1}$ is a graded identity for $R=A\otimes B$. 

($\Leftarrow$)

For each $h\in H$ we choose an element $b_h$  in $B_h$. There exists a substitution $x_i\mapsto a_i\in A_e$ in $f\mathrm{Cap}_{t-1}$ that results in a non-zero element in $A_e$. Denote $a$ the result of this substitution in $\mathrm{Cap}_{t-1}$. Then $ca\neq 0$, where $c=f(a_1,\dots, a_m)$.  We now obtain substitutions in each $f^{h_{i}}$ that result in the element $c\otimes b_{h_i}$, for $i=1,\dots, u$. Hence we obtain a substitution in  $g(f^{h_1},\dots, f^{h_{u}})\mathrm{Cap}_{t-1}$ that results in $$g(c\otimes b_{h_1}, \cdots, c\otimes b_{h_u})(a\otimes 1)=c^ua\otimes g(b_{h_1},\cdots, b_{h_u}).$$ Since $f$ is an $A$-good polynomial the element $c$ lies in $Z(A_e)$. Therefore $c$ a diagonal matrix, hence from $ca\neq 0$ we conclude that $c^{u}a\neq 0$. Therefore $g(b_{h_1},\cdots, b_{h_u})=0$. Since the elements $b_{h_1},\dots, b_{h_u}$ are arbitrary we conclude that $g(x_1^{h_1},\dots, x_u^{h_u})$ is a graded identity for $R$.
\hfill $\Box$

\begin{lemma}\label{coid}
	Let $\mathbb{K}$ be a field of characteristic zero and let $A=\oplus_{g\in G}A_g$ be a $\mathbb{K}$-algebra with a $G$-grading. If the algebras $B$ and $B^{\prime}$ satisfy the same ordinary identities then $A\otimes B$  and $A\otimes B^{\prime}$ satisfy the same graded identities.
\end{lemma}
\textit{Proof.}
Let $F=\mathbb{K}\langle X \rangle/Id(B)$, we prove first that $Id_G(A\otimes B)=Id_G(A\otimes F)$. Let $f(x_1^{g_1},\cdots, x_n^{g_n})$ be a multilinear graded identity for $A\otimes F$ and let $a_1\otimes b_1,\dots, a_n\otimes b_n$ be elements of $A\otimes B$ such that $a_i\otimes b_i\in (A\otimes B)_{g_i}$, $i=1,\dots, n$. Let $\varphi:A\otimes F\rightarrow A\otimes B$ be a graded homomorphism such that $\varphi(a_i\otimes \overline{x_i})=a_i\otimes b_i$. We have $f(a_1\otimes \overline{x_1},\dots, a_n\otimes \overline{x_n})=0$, therefore \[0=\varphi(f(a_1\otimes \overline{x_1},\dots, a_n\otimes \overline{x_n}))=f(a_1\otimes b_1,\dots, a_n\otimes b_n).\] Hence we conclude that $f$ is a graded identity for $A\otimes B$. Since $\mathrm{char}\ \mathbb{K}=0$ this implies that $Id_G(A\otimes F)\subseteq Id_G(A\otimes B)$. Now let $f(x_1^{g_1},\cdots, x_n^{g_n})$ be a multilinear graded identity for $A\otimes B$ and let $\varphi:\mathbb{K}\langle X_G \rangle \rightarrow A\otimes F$ be a graded homomorphism. We write\[\varphi(f)=a_1\otimes f_1+\cdots + a_m\otimes f_m,\] where $a_1,\dots, a_m$ are linearly independent elements of $A$. If $f_i\neq 0$ for some $i$ then there exists a homomorphism $\psi_0:F\rightarrow B$ such that $\psi_0(f_i)\neq 0$. Then $\psi=Id_A\otimes \psi_0$ is a graded homomorphism from $A\otimes F$ to $A\otimes B$ and \[\psi\circ\varphi(f)=a_1\otimes \psi_0(f_1)+\cdots + a_m\otimes \psi_0(f_m)\neq 0,\] this is a contradiction because  $\psi\circ\varphi:\mathbb{K}\langle X_G \rangle \rightarrow A\otimes B$ is a graded homomorphism and $f$ is a graded polynomial identity for $A\otimes B$. Therefore $f_1=f_2=\dots=f_m=0$. Hence $f$ is a graded identity for $A\otimes F$. This proves that $Id_G(A\otimes B)\subseteq Id_G(A\otimes F)$. Therefore  $Id_G(A\otimes B)=Id_G(A\otimes F)$. Clearly we also have  $Id_G(A\otimes B^{\prime})=Id_G(A\otimes F)$.
\hfill $\Box$

\begin{lemma}\label{elemc}
	Let $R=A\otimes B$, where $A=UT(d_1,\dots, d_m)$ has the elementary grading induced by the $n$-tuple $(g_1,\dots, g_n)\in G^n$ and $B$ is an algebra that satisfies the same ordinary identities as $M_t(\mathbb{K})$ and has a $G$-grading with support the subgroup $H$ of $G$. Let $\alpha:G\rightarrow G/H$ denote the canonical epimorphism. The coarsening $^{\alpha}R$ satisfies the same $G/H$-graded identities as the algebra $UT(d_1t,\dots, d_mt)$ with the elementary grading induced by $(g_1H, g_1H,\dots, g_nH, g_nH)$, where each lateral class is repeated $t$ times.
\end{lemma}
\textit{Proof.}
Note that $^{\alpha}R=(^{\alpha}A)\otimes B$, where $B$ has the trivial $G/H$-grading. Lemma \ref{coid} implies that $^{\alpha}R$ satisfies the same graded identities as $(^{\alpha}A)\otimes M_t(\mathbb{K})$. Let $\varphi:M_n(\mathbb{K})\otimes M_t(\mathbb{K})\rightarrow M_{nt}(\mathbb{K})$ denote the canonical isomorphism: if we write the matrices in $M_{nt}(\mathbb{K})$ as $n^2$ blocks of size $t\times t$ then given $a=(a_{ij})\in M_n(\mathbb{K})$ and $b\in M_t(\mathbb{K})$ the $(i,j)$-th block of the matrix $\varphi(a\otimes b)$ is $a_{ij}b$. We consider $A=UT(d_1,\dots, d_m)$ as an ungraded subalgebra of $M_n(\mathbb{K})$, where $n=d_1+\cdots +d_m$. The restriction of $\varphi$ to $A\otimes M_t(\mathbb{K})$ is an isomorphism (of ungraded algebras) onto the subalgebra $UT(d_1t,\dots, d_mt)$ of $M_{nt}(\mathbb{K})$. If $UT(d_1t,\dots, d_mt)$ has the elementary grading induced by $(g_1H, g_1H,\dots, g_nH)$, where each lateral class is repeated $t$ times, then this restriction is an isomorphism of graded algebras from $(^{\alpha}A)\otimes M_t(\mathbb{K})$ to $UT(d_1t,\dots, d_mt)$.
\hfill $\Box$

The next two results will be used in the proof of Theorem \ref{main}.

\begin{proposition}\cite[Proposition 1]{BD}\label{propquot}
	Let $A$ and $B$ be two algebras $G$-graded algebras and $\alpha:G\rightarrow H$ a homomorphism of groups. If $Id_G(A)=Id_G(B)$  then $Id_H(^{\alpha}A)=Id_H(^{\alpha}B)$.
\end{proposition}


In order to proceed to our main result we need to describe the ordinary identities of graded division algebras. 

\begin{theorem}\cite[Theorem 3]{DS}
	Let $G$ be a solvable group and suppose that either $\mathbb{K}$ has
	characteristic $0$ or $\mathbb{K}$ has characteristic $p>0$ and $G$ has no elements of order $p$. Then $K^{\alpha}G$ is semisimple.
\end{theorem}

A semisimple algebra is the subdirect product of primitive algebras. It follows from Kaplansky's Theorem that it satisfies the same identities of a matrix algebra $M_t(\mathbb{K})$, for some $t$ (see, for example, the remarks after Theorem 3.2 in \cite{ABK}). This discussion proves the following corollary.

\begin{corollary}\label{iddiv}
	Let $G$ be an abelian group and $\mathbb{K}$ an algebraically closed field of characteristic zero. If $D$ is a graded division algebra that satisfies an ordinary polynomial identity then there exists a natural number $t$ such that $D$ satisfies the same ordinary identities as $M_t(\mathbb{K})$.
\end{corollary}

We recall that the p.i. degree of a p.i. algebra $D$ is the largest integer $d$ such that the multinear identities of $B$ follow from the multilinear identities of $M_d(\mathbb{K})$ (see \cite[pg. 90, Definition 7.1.8]{DF}). Therefore the natural number $t$ in the previous corollary is the p.i. degree of $D$.

\begin{theorem}\label{main}
	Let $G$ be an abelian group. Denote $A=UT(d_1,\dots, d_m)$ and $A^{\prime}=UT(d_1^{\prime},\dots, d_m^{\prime})$ upper block-triangular matrix algebras with elementary gradings induced by the tuples ${\bf g}=(g_1,\dots, g_n)$ and ${\bf g}^{\prime}=(g_1^{\prime},\dots, g_n^{\prime})$, respectively. Let $K$ be a finitely generated subgroup of $G$ and let $B$ and $B^{\prime}$ be $K$-graded division algebras that satisfy an ordinary polynomial identity and let $R=A\otimes B$, $R^{\prime}=A^{\prime}\otimes B^{\prime}$. If $Id_G(R)=Id_G(R^{\prime})$ then $\mathrm{supp}\ B = \mathrm{supp}\ B^{\prime}$ and $B\cong B^{\prime}$. If $B$ and $B^{\prime}$ have p.i. degree $t$ then $Id_{G/H}(^{\alpha}R)=Id_{G/H}(U)=Id_{G/H}(U^{\prime})=Id_{G/H}(^{\alpha}R^{\prime})$, where $\alpha:G\rightarrow G/H$ is the canonical epimorphism, and $U$, $U^{\prime}$ denote the elementary gradings induced by $(g_1H, g_1H,\dots, g_nH, g_nH)$ and $(g_1^{\prime}H, g_1^{\prime}H,\dots, g_n^{\prime}H, g_n^{\prime}H)$, respectively, on the algebra $UT(d_1t,\dots, d_mt)$, where each lateral class is repeated $t$ times. Moreover, if $U\cong U^{\prime}$ then $R\cong R^{\prime}$.
\end{theorem}
\textit{Proof.}
Corollary \ref{supp} and Proposition \ref{div} imply that $\mathrm{supp}\ B= \mathrm{supp}\ B^{\prime}$ and $Id_H(B)=Id_H(B^{\prime})$, therefore Remark \ref{idisgdiv} implies that $B\cong B^{\prime}$. Corollary \ref{iddiv} and Lemma \ref{elemc} imply that $^{\alpha}R$, $^{\alpha}R^{\prime}$ satisfy the same $G/H$-graded identities as $U$ and $U^{\prime}$ respectively. The equality $Id_{G/H}(^{\alpha}R)=Id_{G/H}(^{\alpha}R^{\prime})$ follows from Proposition \ref{propquot}. Therefore $Id_{G/H}(^{\alpha}R)=Id_{G/H}(U)=Id_{G/H}(U^{\prime})=Id_{G/H}(^{\alpha}R^{\prime})$. Now we assume that $U$ and $U^{\prime}$ are isomorphic as $G/H$-graded algebras. Corollary \ref{appl} imply that there exist a $gH\in G/H$,  $\sigma\in S_{td_1}\times \dots \times S_{td_m}$ such that $g_i^{\prime}H=g_{\sigma(i)}H(gH)$ for $i=1,\dots,tn$. Therefore there exists $\tau\in S_{d_1}\times \dots \times S_{d_m}$ such that $g_i^{\prime}H=g_{\tau(i)}H(gH)$, for $i=1,\dots, n$. Clearly this implies that that there exists $h_1,\dots, h_n \in H$ such that $g_i^{\prime}=g_{\tau(i)}h_{\tau(i)}g$, for $i=1,\dots,n$. Hence it follows from Corollary \ref{appl} that $R\cong R^{\prime}$. 
\hfill $\Box$

In \cite{DS} the authors provide sufficient conditions to ensure that the algebras $U$ and $U^{\prime}$ in the previous theorem are isomorphic as graded algebras. To state their result we need the notion of invariance subgroup of an elementary grading introduced in that paper.

\begin{definition}
	Let $A$ denote the matrix algebra $M_n(\mathbb{K})$ with the elementary grading induced by $(g_1,\dots, g_n)$. The map $\omega_A:G\rightarrow \mathbb{N}$ associates to each $g\in G$ the number of indices $i$ such that $g_i=g$. The invariance subgroup of this elementary grading is $$H_A=\{h\in G \mid \omega_A(hg)=\omega_A(g),\  \forall g \in G\}.$$
\end{definition}

\begin{proposition}\cite[Corollary 3.4]{DS}\label{ds}
	 Let $G$ be an abelian group and let $A$ and $B$ be upper block-triangular matrix algebras endowed with an elementary $G$-grading such that $A$ has either at most two block components or finitely many block components but at least one of them, let us say $A_d$, with invariance subgroup $H_A^{(d)}=\langle 1_G \rangle$. Then $A$ and $B$ are isomorphic as $G$-graded algebras if, and only if, $Id_G(A)=Id_G(B)$.
\end{proposition}

The previous proposition and Theorem \ref{main} imply the following result.

\begin{corollary}
Let $G$ be an abelian group. Denote $A=UT(d_1,\dots, d_m)$ and $A^{\prime}=UT(d_1^{\prime},\dots, d_m^{\prime})$ upper block-triangular matrix algebras with elementary gradings induced by the tuples ${\bf g}\in G^{n}$ and ${\bf g}^{\prime}$, respectively. Let $K$ be a finitely generated subgroup of $G$, let $B$ and $B^{\prime}$ be $G$-graded division algebras that satisfy an ordinary p.i. and have the same p.i. degree. Let $R=A\otimes B$, $R^{\prime}=A^{\prime}\otimes B^{\prime}$. The algebra $A$ has either at most two block components or finitely many block components but at least one of them has trivial invariance subgroup for $^{\alpha}A$, where $\alpha:G\rightarrow G/H$ is the canonical epimorphism. Then $R$ and $R^{\prime}$ are isomorphic as $G$-graded algebras if, and only if, $Id_G(R)=Id_G(R^{\prime})$.
\end{corollary}
\textit{Proof.}
Let $U$ and $U^{\prime}$ be the elementary graded algebras in the previous theorem. The invariance subgroup of the $i$-th block of $^{\alpha}A$ coincides with the invariance subgroup of the $i$-th block $U$. The number of blocks in $^{\alpha}A$ also coincide with the number of blocks in $U$. Therefore Proposition \ref{ds} implies that $U\cong U^{\prime}$. The result now follows from Theorem \ref{main}.
\hfill $\Box$

\end{document}